\documentclass[12pt]{amsart}
\usepackage[margin=1.0in]{geometry}
\usepackage{amssymb,latexsym, amsthm, mathrsfs, enumerate, fullpage,graphicx}
\usepackage[breaklinks]{hyperref}

\newtheorem{theorem}{Theorem}[section]
\newtheorem{lemma}{Lemma}[theorem]
\newtheorem{rem}{Remark}
\newtheorem*{pf}{Proof}

\title{\large \textbf{ON THE EIGENVALUES OF THE LAPLACIAN ON ELLIPSOIDS OBTAINED AS PERTURBATION OF UNIT SPHERE}}
\author{\small ANANDATEERTHA MANGASULI AND ADITYA TIWARI}
\date{}

\begin{document}
\maketitle
    \textbf{Abstract}. We study the eigenvalues of the Laplacian on ellipsoids that are obtained as perturbations of the standard Euclidean unit sphere in dimension two. A comparison of these eigenvalues with those of the standard Euclidean unit sphere is obtained under a Gaussian curvature condition, in line with the Lichnerowicz theorem on the first positive eigenvalue on a compact Riemannian manifold.

\section{Introduction} \label{intro}
\noindent Let $(N,g)$ be a compact $n$-Riemannian manifold without boundary. Then the Laplacian has the following local expression
\[-\Delta_g: = -\frac{1}{\sqrt{|det(g)|}} {\sum_{i,j=1}^{n}}\frac{\partial}{\partial x^i} \bigg{(}g^{ij} \sqrt{|det(g)|} \frac{\partial}{\partial x^j}\bigg{)},\]
in local coordinates $(U,x^1,x^2,\dots,x^n)$ where $(g^{ij})$ is the inverse of the Riemannian metric $g = (g_{ij})$ in this coordinate system and $det(g)$ is the determinant of $(g_{ij})$. It is well known that the set of the eigenvalues of the Laplacian form a non-decreasing sequence 
\[0 = \Lambda_0 < \Lambda_1\leq \Lambda_2 \leq \Lambda_3 \leq \cdots \cdots \cdots \leq \Lambda_n \leq \cdots \nearrow +\infty,\]
where each eigenvalue has finite multiplicity and is repeated as many times as its multiplicity in the sequence. The eigenspaces are orthogonal to each other and their direct sum is $L^2(N)$. A very famous theorem of Lichnerowicz states that
\begin{theorem}[Lichnerowicz]
Let $(N,g)$ be a compact Riemannian $n$-manifold. If there exists a positive number $k$ such that 
\[Ric_{g} \geq k g,\]
where $Ric_g$ being the Ricci tensor of $g$, then
\[\Lambda_1(g) \geq \frac{n}{n-1} k,\]
where $\Lambda_1(g)$ is the first non-zero eigenvalue of the Laplacian acting on smooth functions on $N$.
\end{theorem}
The eigenvalues of the Laplacian on $S^n$ are explicitly known with the first eigenvalue being $n$. Therefore, the Lichnerowicz theorem compares the first eigenvalue of the Laplacian on any compact $n$-manifold with the first eigenvalue of the Laplacian on the standard unit sphere $(S^n,g_0)$. For compact manifold $(\mathcal{M},g)$ of dimension $n$ with \[Ric_g \geq (n-1)g,\] Lichnerowicz's theorem implies $\Lambda_1(g) \geq \Lambda_1(g_0).$ A natural question which can be asked is ``Does the inequality $\Lambda_i(g) \geq \Lambda_i(g_0)$ holds for all $i \geq 0$?" in the above case.\\
\indent The answer to this question in general is NO. Donnelly\cite{Donnely} constructed a counterexample in dimension four where the above eigenvalue inequality fails. Furthermore, the methods used in constructing the above example can be adopted in higher dimensions but not in dimensions two and three, and hence it became interesting to know about the inequalities of the higher eigenvalues in dimensions two and three. In \cite{ATM2}, Mangasuli showed that for all $j\geq 1$, the above inequality holds for left-invariant metrics on $S^3$. Furthermore, in \cite{ATM1}, he showed that if the underlying manifold is $S^2$, then this inequality is also satisfied for certain metrics sufficiently close to $g_0$, which are obtained by analytically perturbing $g_0$ through rotationally symmetric and conformal metrics. \\
\indent In this article, we show given $L\geq 1$, there exist $j\in \mathbb{N}$ such that the above inequality holds for all $i \leq j$, when the underlying manifold is an ellipsoid with the Gaussian curvature greater than or equal to one at every point.\\
The first section contains a motivation behind the problem, while the notation and preliminaries are discussed in section two. The article concludes with stating and proving the main result in section three.

\section{Ellipsoid}    
\label{Gcurv}
In case of $2$-manifolds, $Ric_g$ equals to $K g$, where $K$ denotes the Gaussian curvature. Therefore, $2$-manifolds with $K\geq 1$, satisfy $Ric_g \geq (2-1) g$. In this section, we define the notations and obtain certain conditions on semi-axes lengths for the Gaussian curvature of the ellipsoid to be greater than equal to one at all points. An ellipsoid is a 2-dimensional manifold represented as: 
\[\bigg{\{}(x,y,z) \in \mathbb{R}^3 \bigg{|} \frac{x^2}{A^2} + \frac{y^2}{B^2} + \frac{z^2}{C^2} =1 \bigg{\}},\]
where $A$, $B$ and $C$ $>$ 0 are semi-axes lengths of the ellipsoid. We will use the following notation to represent an ellipsoid. Given $A, B$ and  $C >0$,
\[E_{A, B, C}:=\bigg{\{}(x,y,z) \in \mathbb{R}^3 \bigg{|} \frac{x^2}{A^2} + \frac{y^2}{B^2} + \frac{z^2}{C^2} =1 \bigg{\}}.\]
When  all the semi-axes lengths are unequal, the corresponding ellipsoid is called a \textit{triaxial ellipsoid}. 
When any of the two semi-axes lengths are equal, we will call the ellipsoid a \textit{biaxial ellipsoid}. When all are equal, it is a \textit{sphere}.
Let $p = (x,y,z) \in E_{A,B,C}$, then the Gaussian curvature at $p$ is given by
\begin{equation}
    K(p)=\frac{1}{A^2B^2C^2(\frac{x^2}{A^4}+\frac{y^2}{B^4}+\frac{z^2}{C^4})^2},
\end{equation}
where $A , B$ and $C > 0.$ We will find the minimum of the Gaussian curvature on the ellipsoid and use it to compute the condition on axes lengths for $K \geq 1$.
\subsection{Maxima and Minima of Curvature on $E_{A,B,C}$:}  
Using the method of Lagrange's multipliers to find the minimum curvature on the ellipsoid $E_{A,B,C}$, we get the points of extrema as 
\[\{(\pm A,0,0) , (0,\pm B , 0) , (0 , 0, \pm C )\},\]
with corresponding extrema values as
\[h(\pm A,0,0) = \frac{A^2}{B^2C^2},\]
\[h(0,\pm B,0) = \frac{B^2}{A^2C^2},\]
\[h(0,0,\pm C) = \frac{C^2}{A^2B^2}.\]
In order to obtain the minimum curvature, we specify an ordering of the axes lengths $A$, $B$ and $C$. Suppose $A \geq B \geq C >0$. Then,
\begin{eqnarray*}
K_{max}(E_{A,B,C}):=\max_{p\in E_{A,B,C}} K(p) = K(\pm A,0,0) =  \max_{(x,y,z) \in S}h(x,y,z) = \frac{A^2}{B^2C^2},\\
K_{min}(E_{A,B,C}):=\min_{p \in E_{A,B,C}} K(p) = K(0,0,\pm C) = \min_{(x,y,z) \in S}h(x,y,z) = \frac{C^2}{A^2B^2}.
\end{eqnarray*} 
We revise our notation as follows: $E_{A,B,C}$ will denote all the ellipsoids with semi-axes lengths $A\geq B \geq C> 0$ i.e; $E_{A,B,C}$ will represent six different ellipsoids. The minimum and maximum curvatures of these ellipsoids are given as, $K_{min}(E_{A,B,C})=\frac{C^2}{A^2B^2}$ and $K_{max}(E_{A,B,C})= \frac{A^2}{B^2C^2}$. With the revised notation, there will be two types of biaxial ellipsoids $E_{A,B,B}$ and $E_{A,A,B}$, where $A \geq B >0$.

\subsection{Condition on semi-axes lengths for $K_{min} \geq 1$.}
\subsubsection{Biaxial ellipsoids of type $E_{A,B,B}$.}
We have $K_{min}(E_{A,B,B}) = \frac{1}{A^2}.$
Since $A > B >0$, therefore
\begin{equation} \label{A}
K_{min}(E_{A,B,B}) \geq 1 \iff 1\geq A > B>0.    
\end{equation}
\subsubsection{Biaxial ellipsoids of type $E_{A,A,B}$.}
We have $K_{min}(E_{A,A,B}) = \frac{B^2}{A^4}.$ Therefore,
\begin{equation} \label{B}
    K_{min}(E_{A,A,B}) \geq 1 \iff 1\geq A> B \geq A^2 >0.
\end{equation} 
\subsubsection{Triaxial ellipsoids of type $E_{A,B,C}$.}
We have $K_{min}(E_{A,B,C}) = \frac{C^2}{A^2 B^2}.$ Therefore,
\begin{equation} \label{C}
    K_{min}(E_{A,B,C}) \geq 1 \iff 1\geq A > B > C \geq A B >0.
\end{equation}

\subsection{Ellipsoids as perturbations of Sphere} \label{sec4}
Let $A:=1+\alpha \epsilon \geq $ $B:=1+\beta \epsilon \geq $ $C:= 1 + \gamma \epsilon > 0 $ for $\epsilon>0$ where $\alpha, \beta $ and $\gamma$ are fixed real numbers. Then the family of ellipsoids represented as below
\[E_{A,B,C}(\epsilon) := \bigg{\{} (x,y,z) \in \mathbb{R}^3 \bigg{|}\frac{x^2}{(1+\alpha \epsilon)^2} + \frac{y^2}{(1+\beta \epsilon)^2} + \frac{z^2}{(1+\gamma \epsilon)^2}=1\bigg{\}},\]
is close to $E_{1,1,1}$ for $\epsilon$'s sufficiently close $0$ with $E_{A,B,C}(0) = E_{1,1,1}$. For $\epsilon>0$, inequalities $(\ref{A})$, $(\ref{B})$ and $(\ref{C})$ from the above section implies that 
\begin{equation} \label{ABC}
     0 \geq \alpha >  \beta  > \gamma  >-1/\epsilon.
\end{equation}
Using inequalities derived in the above section and some calculus, we will compute sufficient condition on $\alpha, \beta,$ $\gamma$ and $\epsilon$ for the Gaussian curvature to be greater than equal to one at all points in $E_{A,B,C}(\epsilon)$. We will write $K_{min}(E_{A,B,C})=:K_{min}(\epsilon)$.

\subsubsection{Biaxial ellipsoids of type $E_{A,B,B}(\epsilon)$.}
Inequality $(\ref{A})$ and $(\ref{ABC})$ implies $0 \geq \alpha > \beta > -1/\epsilon$.
Since $K_{min}(\epsilon) = \frac{1}{(1+\alpha \epsilon)^2}$ and $\epsilon > 0$,
 therefore
\begin{equation} \label{Bi1}
K_{min}(\epsilon) \geq 1 \; \text{and} \; A >  B > 0 \iff \quad \begin{cases}
           \alpha = 0 > \beta\;\text{and}\;0<\epsilon<-1/\beta,\\
           0 > \alpha > \beta\; \text{and} \; \alpha/\beta < \epsilon < -1/\beta.
             \end{cases}   
\end{equation}

\subsubsection{Biaxial ellipsoids of type $E_{A,A,B}(\epsilon)$.}
Since $K_{min}(\epsilon) = \frac{(1+\epsilon \beta)^2}{(1+ \epsilon \alpha)^4}$,
therefore
\begin{equation} \label{Bi2}
K_{min}(\epsilon) \geq 1 \; \text{and} \; A >  B > 0 \iff \quad \alpha < 0\;\text{and}\; 2 \alpha < \beta < \alpha\;\text{and}\;0<\epsilon\leq \frac{\beta-2 \alpha}{\alpha^2}.    
\end{equation}

\subsubsection{Triaxial ellipsoids }
$K_{min}(\epsilon) =\frac{(1+\gamma \epsilon)^2}{(1+\beta \epsilon)^2 (1 + \gamma \epsilon)^2 }$. Therefore,
\begin{equation} \label{Triellip}
K_{min}(\epsilon) \geq 1\; \text{and}\; A > B > C > 0 \iff 0 > \alpha > \beta > \gamma > \alpha + \beta \; \text{and} \; 0 < \epsilon \leq \frac{\gamma - \alpha -\beta}{\alpha \beta},    
\end{equation}
\noindent where $A = 1 + \alpha \epsilon$, $B = 1 + \beta \epsilon$ and $C = 1 + \gamma \epsilon$.

\section{Comparison of the Eigenvalues}
In this section, we will state and prove the main theorem using results from \cite{es2021}.
\subsection{Biaxial Ellipsoids}
Assume $A , B > 0$ and $E_{A,A,B}$ be the biaxial ellipsoid with no-ordering in $A$ and $B$. Also, assume that $A=1 + \alpha \epsilon$ and $B = 1+ \beta \epsilon$ and $-\Delta_g$ is the positive Laplace Beltrami operator on $E_{A,A,B}$. Let $g_{\epsilon}$ be the Riemannian metric obtained by restricting the Euclidean metric from $\mathbb{R}^3$ to $E_{A,A,B}$. With these notations, Theorem 2 of \cite{es2021} is restated as follows.
\begin{theorem} \label{Perturbegvalue}
Given $L \in \mathbb{N}$ and $\alpha$, $\beta \in \mathbb{R}$. There exists $\epsilon_0(\alpha,\beta,L)$ such that for all $0<\epsilon < \epsilon_0$ and $\Lambda(\epsilon) \in spec(-\Delta_{g_{\epsilon}})\cap [0,L(L+1)]$, we have
\begin{equation*}
    \Lambda(\epsilon) = l(l+1) + \lambda_1 \epsilon + O(\epsilon^2), 
\end{equation*}
for $l = 0,1,2,\dots,L$ and $m =-l,\dots,l$ with $\lambda_1$ given by the formula
\begin{equation*} \label{evalue1}
    \lambda_1 = -2 \alpha l (l+1) + (\alpha-\beta) \frac{2 l( l +1)}{(2 l+3)(2 l-1)}(2 l^2-2m^2+2 l-1).
\end{equation*}
Moreover, each $\Lambda(\epsilon)$ has multiplicity two except for those where $m=0$, which in this case corresponds to multiplicity one.
\end{theorem}
We define a new notation for the eigenvalues of the ellipsoid $E_{A,A,B}(\epsilon)$ with above theorem in mind. For $L\in \mathbb{N}$ and $\alpha, \beta\in \mathbb{R}$, write $\Lambda_l^m(\epsilon)$ for each eigenvalue corresponding to $l \in \{0,1,\dots,L\}$ and $m \in \{-l,\dots,0,\dots,l\}$ and the corresponding first order term $\lambda_1$ by $\lambda_1(l,m)$.
\begin{rem}
In Theorem $(\ref{Perturbegvalue})$, the semi-axes lengths are not ordered. Therefore, for ellipsoids of the type $E_{A,B,B}$, by interchanging $\alpha$ and $\beta$, $\lambda_1(l,m)$ is given by
\begin{equation} \label{evalue2}
    \lambda_1(l,m) = -2 \beta l (l+1) + (\beta-\alpha) \frac{2 l( l +1)}{(2 l+3)(2 l-1)}(2 l^2-2m^2+2 l-1).
\end{equation} 
\end{rem}
We prove the following lemmas to be used later in proving the main theorem.
\begin{lemma} \label{Eaab}
Let $E_{A,A,B}(\epsilon)$ be a family of biaxial ellipsoids with the Gaussian curvature greater than equal to one at all points. Given any $L \in \mathbb{N}$, there exists an $\epsilon_1(A,B,L)>0$, such that for all $0<\epsilon<\epsilon_1$
\[\lambda_1(l,m) \geq 0, \qquad \text{for}\;\;\text{all} \;\;l=0,1,\dots,L, \qquad \text{and} \quad m -l,\dots,l.\]
\begin{pf}
Using $(\ref{Bi1})$, we have $ 0 \geq \alpha > \beta $ as $K_{min}(\epsilon) \geq 1$. Therefore
\[\lambda_1(l,m) = -2 \alpha l (l+1) + (\alpha-\beta) \frac{2 l( l +1)}{(2 l+3)(2 l-1)}(2 l^2-2m^2+2 l-1) \geq 0,\]
for all $l = 0,1,\dots,L$ and $m = -l,\dots,0,\dots,l$.
\end{pf}
\end{lemma}

\begin{lemma} \label{Eabb}
Let $E_{A,B,B}(\epsilon)$ be a family of biaxial ellipsoids with the Gaussian curvature greater than equal to one at all points. Given any $L \in \mathbb{N}$, there exists an $\epsilon_1(A,B,L)>0,$ such that for all $\epsilon<\epsilon_1$
\[\lambda_1(l,m) \geq 0, \qquad \text{for}\;\;\text{all} \;\;l=0,1,\dots,L, \qquad \text{and} \quad m = -l,\dots,l.\]
\begin{pf}
We have $0 > \alpha > \beta $ by $(\ref{Bi2})$ as $K_{min}(\epsilon) \geq 1$. Writing $\beta = \alpha -\delta $ where $\delta > 0$. Using $(\ref{evalue2})$ and replacing $\beta $ with $\alpha -\delta$, we have
\begin{equation*}
    \lambda_1(l,m) = -2 \alpha l (l+1) + 2 \delta l (l+1)  \frac{(2 l^2 + 2 l - 2 + 2 m^2)}{(2l+3)(2l-1)}\geq 0.
\end{equation*}

\end{pf}
\end{lemma}
\subsection{Triaxial ellipsoids}
We will restate Theorem $4$ of \cite{es2021}. Assume $\epsilon > 0$, $A = 1 + \alpha \epsilon$, $B = 1 + \beta \epsilon$ and $C = 1 + \gamma \epsilon$ such that $A , B$ and $C > 0$.
\begin{theorem} \label{Trispec}
Given $L \in \mathbb{N}$ and $\alpha$, $\beta$, $\gamma \in \mathbb{R}$ with atleast one being non-zero and $g_{\epsilon}$  the metric from $\mathbb{R}^3$ restricted to $E_{A,B,C}$. Then there exist $\epsilon_0(\alpha,\beta,\gamma,L)$ such that for all $0<\epsilon < \epsilon_0$ and $\Lambda(\epsilon) \in$ spec($-\Delta_g$)$\cap [l(l+1)-2l , l(l+1)+2l)]$, we have 
\begin{equation*} \label{egtri}
    \Lambda(\epsilon) = l(l+1)+ \lambda_1 \epsilon + \mathcal{O}(\epsilon^2),
\end{equation*}
where $\lambda_1$ is an eigenvalue of $2 l +1$ eigenvalues of the following symmetric tridiagonal matrices denoted by $M_{\cos,e}$, $M_{\cos,o}$, $M_{\sin,e}$ and $M_{\sin,o}$.
\end{theorem}
We need some notations to define these matrices. Define
\[l_e= \begin{cases}
l & ,\text{if $l$ is even},\\
l-1 & ,\text{if $l$ is odd}.
           \end{cases}; \qquad l_o= \begin{cases}
l &, \text{if $l$ is odd},\\
l-1 &, \text{if $l$ is even}.
           \end{cases}\]
Then $M_{\cos,e}$ is given by
\[M_{\cos,e} = \begin{bmatrix}
a_0 && b_0 && 0     &&\dots && 0\\
b_0 && a_1 && b_1 &&      && 0\\
0  && b_1 && a_2 &&      && \vdots  \\
\vdots &&  &&  &&\ddots &&    \vdots\\
 && && && && b_{k-1}\\
0 && \dots && && b_{k-1} && a_k
\end{bmatrix}_{(k+1)\times (k+1)}, \]
for $l =2k$ or $2k+1$ i.e; $l_e = 2k$ and
 $M_{\sin,e}$ is given as
\[M_{\sin,e} = \begin{bmatrix}
a_1 && b_1 && 0     &&\dots && 0\\
b_1 && a_2 && b_2 &&      && 0\\
0  && b_2 && a_3 &&      && \vdots  \\
\vdots &&  &&  &&\ddots &&    \vdots\\
 && && && && b_{k-1}\\
0 && \dots && && b_{k-1} && a_k
\end{bmatrix}_{k\times k}, \]
for $l =2k$ or $2k+1$ i.e; $l_e = 2k$ and
\begin{align} \label{Me}
-a_p &= 2 \gamma l (l+1) + (\alpha + \beta - 2 \gamma)\frac{2 l (l+1)}{(2l-1)(2l+3)}(l^2+4 p^2 + l-1),\quad p =0 ,1,\dots,k,\\
-b_0 &= \sqrt{2} (\beta-\alpha)\frac{ l (l+1)}{(2 l +3) (2 l -1)}\sqrt{(l-1) l (l+1) (l+2)},\\
-b_p &= (\beta-\alpha)\frac{ l (l+1)}{(2 l +3) (2 l -1)}\sqrt{(l-2 p -1) (l-2p) (l+2p+1) (l+2p+2)},\quad p = 1,\dots,k-1.
\end{align}
When $l=2k+1$, $M_{\cos,o}$ is given by
\[M_{\cos,o} = \begin{bmatrix}
{\phi}_0 && {\psi}_0 && 0     &&\dots && 0\\
{\psi}_0 && {\phi}_1 && {\psi}_1 &&      && 0\\
0  && {\psi}_1 && {\phi}_2 &&      && \vdots  \\
\vdots &&  &&  &&\ddots &&    \vdots\\
 && && && && {\psi}_{k-1}\\
0 && \dots && && {\psi}_{k-1} && {\phi}_k
\end{bmatrix}_{(k+1)\times (k+1)}, \]
and when $l=2k$, $M_{\cos,o}$ is given by 
\[M_{\cos,o} = \begin{bmatrix}
{\phi}_0 && {\psi}_0 && 0     &&\dots && 0\\
{\psi}_0 && {\phi}_1 && {\psi}_1 &&      && 0\\
0  && {\psi}_1 && {\phi}_2 &&      && \vdots  \\
\vdots &&  &&  &&\ddots &&    \vdots\\
 && && && && {\psi}_{k-2}\\
0 && \dots && && {\psi}_{k-2} && {\phi}_{k-1}
\end{bmatrix}_{k\times k} .\]
$M_{\sin,o}$ is given by
\[M_{\sin,o} = \begin{bmatrix}
\Tilde{{\phi}_0} && {\psi}_0 && 0     &&\dots && 0\\
{\psi}_0 && {\phi}_1 && {\psi}_1 &&      && 0\\
0  && {\psi}_1 && {\phi}_2 &&      && \vdots  \\
\vdots &&  &&  &&\ddots &&    \vdots\\
 && && && && {\psi}_{k-1}\\
0 && \dots && && {\psi}_{k-1} && {\phi}_k
\end{bmatrix}_{(k+1)\times (k+1)}, \]
when $l=2k+1$ and
\[M_{\sin,o} = \begin{bmatrix}
\Tilde{{\phi}_0} && {\psi}_0 && 0     &&\dots && 0\\
{\psi}_0 && {\phi}_1 && {\psi}_1 &&      && 0\\
0  && {\psi}_1 && {\phi}_2 &&      && \vdots  \\
\vdots &&  &&  &&\ddots &&    \vdots\\
 && && && && {\psi}_{k-2}\\
0 && \dots && && {\psi}_{k-2} && {\phi}_{k-1}
\end{bmatrix}_{k\times k}, \]
when $l = 2k$, where
\begin{align} \label{phip1}
-\phi_0 &= (\frac{3 \alpha}{2} + \frac{\beta}{2} - 2\gamma) \frac{2 l^2 (l+1)^2}{(2 l -1)(2 l +3)} + 2 l (l+1) \gamma,\\
-\Tilde{\phi_0} &= (\frac{3 \beta}{2} + \frac{\alpha}{2} - 2\gamma) \frac{2 l^2 (l+1)^2}{(2 l -1)(2 l +3)} + 2 l (l+1) \gamma,
\end{align}
\begin{align*}
    -\phi_p = (\alpha + \beta) \frac{2 l (l+1)}{(2l+3)(2l-1)}(l^2 + 4 p^2 + 4 p +l)
\end{align*}
\begin{equation*}
    + 2 \gamma \frac{l (l+1)}{(2 l-1)(2 l+3)}(2 l^2 - 8 p^2 - 8p+2l -3),\qquad p=1,\dots,k,
\end{equation*}
\begin{align*}
-\psi_p &=\frac{(\beta-\alpha) l (l+1)}{(2 l +3) (2 l -1)}\sqrt{(l-2 p -2) (l-2p-1) (l+2p+2) (l+2p+3)}, \quad  \;\;p=1,\dots,k-1.
\end{align*}
\subsubsection{Positive Definiteness of Tridiagonal Matrices.}
A real symmetric matrix is positive definite if determinant of all principal minors is positive. Milica Andelie and C.M. da Fonseca in \cite{PD} has stated the following result for positive definiteness of real symmetric matrices with positive diagonal entries.  
\begin{theorem} \label{PDM}
Let $A_n$ be a real symmetric tridiagonal matrix with positive diagonal entries defined as:
\[A_n := \begin{bmatrix}
r_1 && s_1 && 0     &&\dots && 0\\
s_1 && r_2 && s_2 &&      && 0\\
0  && s_2 && r_3 &&      && \vdots  \\
\vdots &&  &&  &&\ddots &&    \vdots\\
 && && && && s_{n-1}\\
0 && \dots && && s_{n-1} && r_n
\end{bmatrix}_{n\times n}.\]
If 
\begin{equation}
    {s_i}^2 < \frac{1}{4} r_i r_{i+1} \frac{1}{\cos^2{(\frac{\pi}{n+1})}},\quad\text{for}\;\;i=1,\dots,n-1,
\end{equation}
then $A_n$ is positive definite.
\end{theorem}
We will show that the matrices $M_{\cos,e}$, $M_{\cos,o}$, $M_{\sin,e}$ and $M_{\sin,o}$ have positive diagonal entries given the curvature conditions are satisfied and then prove the positive definiteness using the above theorem along with curvature conditions.
\begin{lemma}
Given $\alpha ,\beta$ and $\gamma<0$, $\Tilde{\phi_0}$, $a_p$ and $\phi_p$  are positive for all possible entries of $p$ in above matrices.
\begin{pf}
Rewriting $a_p$ and $\phi_p$ we have,
\begin{equation*}
    -a_p =  (\alpha + \beta)\frac{2 l (l+1)}{(2l-1)(2l+3)}(l^2+4 p^2 + l-1) + 2 \gamma \frac{ l (l+1)}{(2l-1)(2l+3)}(2 l^2-8 p^2 +2 l -1) ,
\end{equation*}
for $p = 0 ,1,\dots,k$, and
\begin{equation*}
    -\phi_p = (\alpha + \beta) \frac{2 l (l+1)}{(2l+3)(2l-1)}(l^2 + 4 p^2 + 4 p +l) + 2 \gamma \frac{l (l+1)}{(2 l-1)(2 l+3)}(2 l^2 - 8 p^2 - 8p+2l -3),
\end{equation*}
for $p = 1,\dots,k.$ Since $l = 2k $ or $2k+1$ and $\alpha , \beta ,\gamma < 0$, therefore
\[a_p > 0, \quad \text{for} \quad p =0 , \dots , k,\]
and
\[\phi_p > 0, \quad \text{for} \quad p =1 , \dots , k.\]
Rewriting $\Tilde{\phi_0}$ and $\phi_0$, we have
\begin{equation*}
    -\Tilde{\phi_0} = (\frac{3 \beta}{2} + \frac{\alpha}{2}) \frac{2 l^2 (l+1)^2}{(2 l -1)(2 l +3)} + 2  \gamma \frac{2 l^2 (l+1)^2}{(2 l -1)(2 l +3)}(2 l^2 + 2 l-3),
\end{equation*}

\begin{equation*}
    -\phi_0 = (\frac{3 \alpha}{2} + \frac{\beta}{2}) \frac{2 l^2 (l+1)^2}{(2 l -1)(2 l +3)} + 2  \gamma \frac{2 l^2 (l+1)^2}{(2 l -1)(2 l +3)}(2 l^2 + 2 l-3).
\end{equation*}
Here $l = 2k $ or $2k+1$ and $\alpha , \beta ,\gamma < 0$. Therefore,
\[\phi_0 > 0 \quad \text{and} \quad \Tilde{\phi_0} > 0.\]
\end{pf}
\end{lemma}
\noindent The order of these matrices depend on $l$ and is given as
\begin{center}
\begin{tabular}{|c|c|c|c|c|}
\hline
$l$ &$M_{\cos,e}$     &  $M_{\sin,e}$ & $M_{\cos,o}$ & $M_{\sin,o}$\\
\hline
$2 k +1$ &$ (k+1) \times (k+1)$     & $k \times k$ & $ (k+1) \times (k+1)$ & $ (k+1) \times (k+1)$\\
\hline
$2 k$ &$ (k+1) \times (k+1)$     & $k \times k$ & $ k \times k$ & $ k \times k$\\
\hline
\end{tabular}    
\end{center}
We will write the order of these matrices for some values of $l$.
\begin{center}
\begin{tabular}{|c|c|c|c|c|}
\hline
$l$ & $M_{\cos,e}$     &  $M_{\sin,e}$ & $M_{\cos,o}$ & $M_{\sin,o}$\\
\hline
$1$ &$1 \times 1$     & - & $1 \times 1$ & $1\times 1$\\
\hline
$2$ &$2 \times 2$     & $1 \times 1$ & $1 \times 1$ & $1\times 1$\\
\hline
$3$ &$2 \times 2$     & $1 \times 1$ & $2 \times 2$ & $2\times 2$\\
\hline
$4$ &$3 \times 3$     & $2 \times 2$ & $2 \times 2$ & $2\times 2$\\
\hline
$5$ &$3 \times 3$     & $2 \times 2$ & $3 \times 3$ & $3\times 3$\\
\hline
$6$ &$4 \times 4$     & $3 \times 3$ & $3 \times 3$ & $3\times 3$\\
\hline
\end{tabular}    
\end{center}
Since all the diagonal entries are positive, matrices of order $1 \times 1$ are positive definite. Therefore, $M_{\cos,e}$ is positive definite for $l=1$, $M_{\sin,e}$ is positive definite for $l = 2$ and $3$, $M_{\cos,o}$ and $M_{\sin,o}$ are positive definite for $l = 1$ and $2$. For matrices of order $2 \times 2$, we will show that the determinant is positive and hence the matrices are positive definite. And for matrices of order greater than $2 \times 2$, we will use Theorem $(\ref{PDM})$ to prove the positive definiteness as diagonal entries are positive.
\begin{lemma} \label{Mcose}
$M_{\cos,e}$ is positive definite for all $l = 1, 2, \dots, L$.
\begin{pf}
Using $(\ref{Me})$, we have
\begin{equation*}
    \frac{4 b_0^2}{ a_0 a_1} = \frac{2 (\beta-\alpha)^2 (l-1) l(l+1) (l+2)}{\{(\alpha+\beta)(l^2+l-1) + \gamma (2 l^2 + 2l -1)\} * \{(\alpha+\beta) (l^2 + l + 3) + \gamma (2 l^2 +2 l -9)\}}.
\end{equation*}
But $0 > \alpha > \beta >\gamma \geq \alpha +\beta$, therefore

\begin{equation*}
\begin{split}
    \frac{4 b_0^2}{ a_0 a_1} &\leq \frac{2 (\beta-\alpha)^2}{\gamma^2} \frac{(l-1) l (l+1) (l+2)}{(3 l^2 + 3 l-2) ( 3 l^2 + 3 l -6)}\\
                            &\leq \frac{2 (\beta-\alpha)^2}{9 \gamma^2} \frac{l (l+1)}{(l-1) (l+2)},\quad \text{for} \quad l\geq 2\\
                            & \leq \frac{4}{9}\frac{ (\beta-\alpha)^2}{\gamma^2}.\\
\end{split}
\end{equation*}
Therefore,
\begin{equation*}
    \frac{4 b_0^2}{ a_0 a_1} \cos^2{(\frac{\pi}{k+2})} \leq \frac{4}{9}\frac{ (\beta-\alpha)^2}{\gamma^2} < 1.
\end{equation*}
For $p > 0$, we have
\begin{equation*}
   \frac{4 b_p^2}{ a_p a_{p+1}} = \frac{ (\beta-\alpha)^2 (l-2 p -1) (l-2p)(l+2p+1) (l+2p+2)}{\{(\alpha+\beta)(l^2+4 p^2 + l-1) + \gamma (2 l^2 - 8 p^2+ 2l -1)\}} 
\end{equation*}
\begin{equation*}
    \times \frac{1}{ \{(\alpha+\beta) (l^2 +4 p^2 +8 p+ l + 3) + \gamma (2 l^2-8 p^2-16 p +2 l -9)\}}.
\end{equation*}
But $0 > \alpha > \beta >\gamma \geq \alpha +\beta$, therefore
\begin{equation*} \label{apap1}
    \begin{split}
    \frac{4 b_p^2}{ a_p a_{p+1}}    & \leq \frac{ (\beta-\alpha)^2}{\gamma^2} \frac{(l+2p+1)(l-2p)}{9 (l-2p-1)(l+2p+2)}\quad \text{for} \quad l \geq 4,\\
     & \leq \frac{2 (\beta-\alpha)^2}{9\gamma^2}.
    \end{split}
\end{equation*}
Therefore, we have
\begin{equation*}
    \frac{4 b_p^2}{ a_p a_{p+1}}  \cos^2{(\frac{\pi}{k+2})} \leq \frac{ 2(\beta-\alpha)^2}{9\gamma^2} < 1\quad \text{where} \quad k\geq 2,
\end{equation*}
for $p = 1 ,\dots, k-1$. Therefore, by Theorem $(\ref{PDM})$, $M_{\cos,e}$ is positive definite for all $l \geq 4$.
When $l=3$ or $2$, $M_{\cos,e}$ will be of order $2 \times 2$. We will show that $a_0 a_1 -b_0^2 >0$ for $l=2$ and $l=3$. Using $(\ref{Me})$, we have 
\begin{eqnarray*}
    a_0 a_1 &=& \frac{4 l^2(l+1)^2}{(2l-1)^2(2l+3)^2}\{(\alpha+\beta)(l^2+l-1) + \gamma (2 l^2 + 2l -1)\}  \{(\alpha+\beta) (l^2 + l + 3) + \gamma (2 l^2 +2 l -9)\}\\
    &>&0.\\
    b_0^2 &=& 2 (\beta-\alpha)^2 \frac{l^2(l+1)^2}{(2l-1)^2(2l+3)^2}(l-1)l(l+1)(l+2) > 0.
\end{eqnarray*}
The determinant $a_0 a_1 - b_0^2$ will be positive if the coefficient of $\alpha^2$ in $a_0 a_1$ is greater than coefficient of $\alpha^2$ in $b_0^2$. Writing the difference of coefficients, we have
\begin{equation*}
     \frac{l^2(l+1)^2}{(2l-1)^2(2l+3)^2} \bigg( 4 (l^2+l-1) (l^2 + l + 3) - (l-1)l(l+1)(l+2)\bigg) > 0 \quad \text{for} \quad l =3\;\text{or}\;4.
\end{equation*}
Therefore, $M_{\cos,e}$ is positive definite for all $l \geq 1$.
\end{pf}
\end{lemma}

\begin{lemma} \label{Msine}
$M_{\sin,e}$ is positive definite for all $l = 2, \dots, L$.
\begin{pf}
By Theorem $(\ref{PDM})$, $M_{\sin,e}$ will be positive definite if
\[4 \frac{b_p^2}{a_p a_{p+1}} \cos^2{(\frac{\pi}{k+1})} < 1 \quad \text{for} \quad p = 1 , \dots, k-1.\]
But $(\ref{apap1})$ implies that
\begin{equation*}
    \frac{4 b_p^2}{ a_p a_{p+1}}
      \leq \frac{ (\beta-\alpha)^2}{\gamma^2} \quad \text{for} \quad l\geq 4.
\end{equation*}
So, we have
\[ 4 \frac{b_p^2}{a_p a_{p+1}} \cos^2{(\frac{\pi}{k+1})} \leq \frac{ (\beta-\alpha)^2}{\gamma^2} < 1. \]
\end{pf}
Therefore, $M_{\sin,e}$ positive definite for all $l \geq 4$. For $l=1$, $M_{\sin,e}$ does not exist and for $l=2$ and $3$, $M_{\sin,e}$ is of order $1 \times 1$ and hence positive definite.
\end{lemma}

\begin{lemma} \label{cosino}
$M_{\cos,o}$ and $M_{\sin,o}$ are positive definite matrices for all $l = 1, 2, \dots, L$.
\begin{pf}
Using $(\ref{phip1})$, we have
\[\frac{4 \psi_0^2}{\phi_0\phi_1} =  \frac{ (\beta-\alpha)^2(l-2) (l-1) (l+2)(l+3)}{\{(3 \alpha+\beta) l (l+1)/2 + \gamma (2 l^2 + 2 l-3)\} \times \{ (\alpha+\beta) (l^2 + +8) + \gamma (2 l^2 + 2 l -19)\}}.\]
But $0 > \alpha > \beta >\gamma \geq \alpha +\beta$, therefore
\begin{equation*}
\begin{split}
    \frac{4 \psi_0^2}{\phi_0\phi_1} 
     & \leq \frac{2 (\beta-\alpha)^2}{15 \gamma^2} \frac{(l-1)(l+2)}{(l-2)(l+3)} \quad \text{for} \quad l\geq 4,\\ 
     & \leq \frac{ (\beta-\alpha)^2}{ \gamma^2}.
\end{split}
\end{equation*}

\noindent Therefore,
\[\frac{4 \psi_0^2}{\phi_0\phi_1} \cos^2{(\frac{\pi}{n+1})} \leq \frac{  (\beta-\alpha)^2}{ \gamma^2} <1,\quad \text{where} \quad n = \begin{cases}
          k+1, & l =2 k+1,\\
          k, & l = 2k.
           \end{cases}\]
Now, for $p>0$,
\begin{equation*} 
    \frac{4 \psi_p^2}{\phi_p \phi_{p+1}} = \frac{ (\beta-\alpha)^2(l-2p-2) (l-2p-1) (l+2p+2)(l+2p+3)}{\{(3 \alpha+\beta) (l^2+4p^2+4p+l) + \gamma (2 l^2 -8p^2-8p+ 2 l-3)\} }
\end{equation*}
\begin{equation*}
    \times \frac{1}{\{ (\alpha+\beta) (l^2 +4p^2+12p+l +8) + \gamma (2 l^2 -8p^2-24p+ 2 l -19)\}}.
\end{equation*}
Therefore,
\begin{equation*}
    \begin{split}
        \frac{4 b_\psi^2}{\phi_p\phi_{p+1}} & \leq \frac{ (\beta-\alpha)^2}{9 \gamma^2} \frac{(l+2p+2)(l-2p-1)}{(l-2p-2)(l+2p+3)}, \quad \text{for} \quad l\geq5,\\
        & \leq \frac{ 2(\beta-\alpha)^2}{9 \gamma^2}.
    \end{split}
\end{equation*}
Hence,
\[\frac{4 \psi_p^2}{\phi_p\phi_{p+1}} \cos^2{(\frac{\pi}{n+1})} \leq \frac{ (\beta-\alpha)^2}{ \gamma^2} < 1, \quad where \quad n = \begin{cases}
          k+1, & l =2 k+1,\\
          k, & l = 2k.
           \end{cases}\]
Therefore, $M_{\cos,o}$ is positive definite for $l \geq 5$ by Theorem $(\ref{PDM})$.
Also, \[\frac{4 \psi_0^2}{\Tilde{\phi_0}\phi_1} =  \frac{ (\beta-\alpha)^2(l-2) (l-1) (l+2)(l+3)}{\{( \alpha+3\beta) l (l+1)/2 + \gamma (2 l^2 + 2 l-3)\} \times \{ (\alpha+\beta) (l^2 + +8) + \gamma (2 l^2 + 2 l -19)\}}.\]
But $0 > \alpha > \beta >\gamma \geq \alpha +\beta$, therefore
\begin{equation*}
\begin{split}
    \frac{4 \psi_0^2}{\Tilde{\phi_0}\phi_1} & \leq \frac{4 (\beta-\alpha)^2}{15 \gamma^2},\quad\text{for} \quad l\geq 4.
\end{split}
\end{equation*}
\noindent Therefore,
\[\frac{4 \psi_0^2}{\Tilde{\phi_0}\phi_1} \cos^2{(\frac{\pi}{n+1})} \leq \frac{ (\beta-\alpha)^2}{ \gamma^2}<1, \quad where \quad n = \begin{cases}
          k+1, & l = 2 k+1,\\
          k, & l = 2k.
           \end{cases}\]
Therefore, $M_{\sin,o}$ is positive definite for all $l\geq 5$ by Theorem $(\ref{PDM})$. When $l=3$ or $4$, we will show that the determinant is positive and hence the matrices are positive definite. The determinant will be positive iff 
\[\frac{3l(l+1)(l^2+l+8)}{2}-(l-2)(l-2)(l+2)(l+3)>0\]
But for $l=3$ and $4$, above quantity is positive. Hence, $M_{\cos,o}$ and $M_{\sin,o}$ are positive definite for $l =3$ and $4$.
\end{pf}
\end{lemma}
Using lemmas (\ref{Mcose}), (\ref{Msine})  and (\ref{cosino}), we have the following theorem. 
\begin{theorem}{} \label{TriM}
$M_{\cos,e}$, $M_{\sin,e}$, $M_{\cos,o}$ and $M_{\sin,o}$ are positive definite for $\alpha$, $\beta$, $\gamma$ and $\epsilon \in (0, \epsilon_1)$ whenever the corresponding triaxial ellipsoids $E_{A,B,C}(\epsilon)$ has the Gaussian curvature greater than equal to one at all points.
\end{theorem}

\begin{theorem}{} \label{Main}
Let $\lambda_i(g_{\epsilon})$ be the $i^{th}$ eigenvalue of the Laplace operator $\Delta_{g_{\epsilon}}$ on the ellipsoid $E_{A,B,C}(\epsilon)$. Consider an analytic perturbation of the standard Euclidean sphere through ellipsoids, $E_{A,B,C}(\epsilon)$, with the Gaussian curvature bounded below by one. Then given $L \in \mathbb{N}$, there exists $j$ depending on $L$ such that 
\[\lambda_{i}(g_{\epsilon}) \geq \lambda_{i}(g_0),\;\;\;\text{for all}\;\;i\leq j,\] for $\epsilon$'s sufficiently close to $0$.
\begin{pf}
The proof follows directly from lemmas $(\ref{Eaab})$, $(\ref{Eabb})$ and Theorem $(\ref{TriM})$.
\end{pf}
\end{theorem}

\bibliography{bi}

\providecommand{\bysame}{\leavevmode\hbox to3em{\hrulefill}\thinspace}
\providecommand{\MR}{\relax\ifhmode\unskip\space\fi MR }
\providecommand{\MRhref}[2]{%
  \href{http://www.ams.org/mathscinet-getitem?mr=#1}{#2}
}
\providecommand{\href}[2]{#2}
\begin{thebibliography}{1}

\bibitem{PD}
Milica An\dj~eli\'{c} and C.~M. da~Fonseca, \emph{Sufficient conditions for
  positive definiteness of tridiagonal matrices revisited}, Positivity
  \textbf{15} (2011), no.~1, 155--159. \MR{2782752}

\bibitem{Donnely}
Harold Donnelly, \emph{Harmonic functions on manifolds of nonnegative {R}icci
  curvature}, Internat. Math. Res. Notices (2001), no.~8, 429--434.
  \MR{1827086}

\bibitem{es2021}
Suresh Eswarathasan and Theodore Kolokolnikov, \emph{Laplace--{B}eltrami
  spectrum of ellipsoids that are close to spheres and analytic perturbation
  theory}, IMA J. Appl. Math. \textbf{87} (2022), no.~1, 20--49. \MR{4370363}

\bibitem{ATM2}
Anandateertha Mangasuli, \emph{On the eigenvalues of the {L}aplacian for
  left-invariant {R}iemannian metrics on {$S^3$}}, Internat. J. Math.
  \textbf{18} (2007), no.~8, 895--901. \MR{2339575}

\bibitem{ATM1}
\bysame, \emph{On the eigenvalues of the {L}aplacian for certain perturbations
  of the standard {E}uclidean metric on {$S^2$}}, Asian J. Math. \textbf{13}
  (2009), no.~2, 271--282. \MR{2559111}

\end{thebibliography}
\bibliographystyle{amsplain}
\end{document}